\documentclass[11pt]{amsart}
\usepackage{amsmath,amsthm,epsfig,amssymb}
\usepackage{xypic}
\input xy
\title{Equilateral triangles in $\mathbb Z^4$  }
\author{Eugen J. Ionascu  }
\curraddr{Department of Mathematics\\ Columbus State University\\4225 University Avenue\\
Columbus, GA 31907} \email{math@ejionascu.edu} \subjclass{}
\date{July  $14^{th}$, 2013}
\begin{document}
\def\sms{\small\scshape}
\baselineskip18pt
%%%%%%%%%%%%%%%%%%%%%%%%%%%%%%%%%%%%%%%%%%%%%%%%%%%%%%%%%%%%%%%%%%%%%%
%%%%%%%%%%%%%%%%%%%% OUR DEFINITIONS %%%%%%%%%%%%%%%%%%%%%%%%%%%%%%%%%
%%%%%%%%%%%%%%%%%%%%%%%%%%%%%%%%%%%%%%%%%%%%%%%%%%%%%%%%%%%%%%%%%%%%%%
\newtheorem{theorem}{\hspace{\parindent}
T{\scriptsize HEOREM}}[section]
\newtheorem{proposition}[theorem]
{\hspace{\parindent }P{\scriptsize ROPOSITION}}
\newtheorem{corollary}[theorem]
{\hspace{\parindent }C{\scriptsize OROLLARY}}
\newtheorem{lemma}[theorem]
{\hspace{\parindent }L{\scriptsize EMMA}}
\newtheorem{definition}[theorem]
{\hspace{\parindent }D{\scriptsize EFINITION}}
\newtheorem{problem}[theorem]
{\hspace{\parindent }P{\scriptsize ROBLEM}}
\newtheorem{conjecture}[theorem]
{\hspace{\parindent }C{\scriptsize ONJECTURE}}
\newtheorem{example}[theorem]
{\hspace{\parindent }E{\scriptsize XAMPLE}}
\newtheorem{remark}[theorem]
{\hspace{\parindent }R{\scriptsize EMARK}}
\renewcommand{\thetheorem}{\arabic{section}.\arabic{theorem}}
\renewcommand{\theenumi}{(\roman{enumi})}
\renewcommand{\labelenumi}{\theenumi}
\newcommand{\Q}{{\mathbb Q}}
\newcommand{\Z}{{\mathbb Z}}
\newcommand{\N}{{\mathbb N}}
\newcommand{\C}{{\mathbb C}}
\newcommand{\R}{{\mathbb R}}
\newcommand{\F}{{\mathbb F}}
\newcommand{\K}{{\mathbb K}}
\newcommand{\D}{{\mathbb D}}
\newcommand{\comments}[1]{}
\def\phi{\varphi}
\def\ra{\rightarrow}
\def\sd{\bigtriangledown}
\def\ac{\mathaccent94}
\def\wi{\sim}
\def\wt{\widetilde}
\def\bb#1{{\Bbb#1}}
\def\bs{\backslash}
\def\cal{\mathcal}
\def\ca#1{{\cal#1}}
\def\Bbb#1{\bf#1}
\def\blacksquare{{\ \vrule height7pt width7pt depth0pt}}
\def\bsq{\blacksquare}
\def\proof{\hspace{\parindent}{P{\scriptsize ROOF}}}
\def\pofthe{P{\scriptsize ROOF OF}
T{\scriptsize HEOREM}\  }
\def\pofle{\hspace{\parindent}P{\scriptsize ROOF OF}
L{\scriptsize EMMA}\  }
\def\pofcor{\hspace{\parindent}P{\scriptsize ROOF OF}
C{\scriptsize ROLLARY}\  }
\def\pofpro{\hspace{\parindent}P{\scriptsize ROOF OF}
P{\scriptsize ROPOSITION}\  }
\def\n{\noindent}
\def\wh{\widehat}
\def\eproof{$\hfill\bsq$\par}
\def\ds{\displaystyle}
\def\du{\overset{\text {\bf .}}{\cup}}
\def\Du{\overset{\text {\bf .}}{\bigcup}}
\def\b{$\blacklozenge$}

\def\eqtr{{\cal E}{\cal T}(\Z) }
\def\eproofi{\bsq}

%%%%%%%%%%%%%%%%%%%%%%%%%%%%%%%%%%%%%%%%%%%%%%%%%%%
%%%%%%%%%%%%%%%%%%%%%%%ABSTRACT%%%%%%%%%%%%%%%%%%%%%%%%%%%%%%%%%%%%
\begin{abstract} We  characterize all three point sets in $\mathbb R^4$ with integer
coordinates, the pairs of which are the same Euclidean distance apart. In three dimensions the problem is characterized in
terms of solutions of the Diophantine equation $a^2+b^2+c^2=3d^2$. In $\mathbb R^4$, our characterization is
essentially based on two different solutions of the same equation. The characterization is existential in nature,
as opposed to the three dimensional situation where we have precise formulae in terms of $a$, $b$ and $c$.
\end{abstract} \maketitle
%%%%%%%%%%%%%%%%%%%%%%%%%%%%%%%%%%%%%%%%%%%%%%%%%%%%%%%%%%%%%%%%%%%%%%%
%%%%%%%%%%%%%%%%%%%%%%%%%%%%%%%%%%%%%%%%%%%%%%%%%%%%%%%%%%%%%%%%%%%%%%%
%%%%%%%%%%%%%%%%%%%%%%%%%%%%%%%%%%%%%%%%%%%%%%%%%%%%%%%%%%%%%%%%%%%%%%%
\section{INTRODUCTION}%%%%%%%%%%%%%%%%%%%%%%%%%%%%%%%%%%%%%%%%%%%%%%%%%

Our work is motivated mainly by the result in  \cite{mc}  and \cite{s} which characterizes the dimensions $n$ in which  $n$ dimensional lattice simplices exist:  {\it all $n$ such that $n+1$ is a sum of 1, 2, 4 or 8 odd squares.}
We are interested in the existence and characterizations of smaller dimension lattice simplices in dimension four. 
The Diophantine equation
\begin{equation}\label{maideophantineeq}
a^2+b^2+c^2=3d^2
\end{equation}

\n   clearly has  trivial solutions such as $a=b=c=d$, but for $d$ odd and greater than one,  it is known (see \cite{ch}, \cite{ejirtlv}) to always have
non-trivial solutions ($(|a|-|b|)^2+(|a|-|c|)^2>0$). An ordered triple of integers $[a, b, c]$
will be called a {\em primitive} solution of (\ref{maideophantineeq}) provided that (\ref{maideophantineeq}) is satisfied,
$\gcd(a,b,c)=1$ and $0<a\le b\le c$. It is easy to see that if $d$ is even then primitive solutions of (\ref{maideophantineeq}) do not exist since
$a$, $b$ and $c$ must be even. We include the primitive solutions corresponding to $d$ odd and $d\le 19$.

$$\begin{tabular}{|c|c|}
  \hline
  % after \\: \hline or \cline{col1-col2} \cline{col3-col4} ...
  d & [a,b,c] \\
   \hline
  1 & [1,1,1]\\
    3 & [1,1,5]\\
      5 & [1,5,7]\\
       7 & [1,5,11]\\
       9 & [1, 11, 11], [5, 7, 13] \\
  \hline
\end{tabular}\ \ \
\begin{tabular}{|c|c|}
  \hline
  % after \\: \hline or \cline{col1-col2} \cline{col3-col4} ...
  d & [a,b,c]\\
   \hline
11& [1, 1, 19], [5, 7, 17], [5, 13, 13]\\
13& [5, 11, 19], [7, 13, 17]\\
15& [1, 7, 25], [5, 11, 23], [5, 17, 19]\\
17 & [1, 5, 29], [7, 17, 23], [11, 11, 25], [13, 13, 23]\\
19& [1, 11, 31], [5, 23, 23], [11, 11, 29], [13, 17, 25]\\
  \hline
\end{tabular}
$$

\vspace{0.1in}

\n As a result of this investigation, we found yet a new  (see \cite{eji}) parametrization of solutions of (\ref{maideophantineeq})
and strong numerical evidence supports the conjecture that it is exhaustive:

{\small \begin{equation}\label{parofmain}
\begin{array}{c} a:=x^2+y^2-z^2-t^2+2(xz-xt-yt-yz),\ \ b:= y^2+t^2-x^2-z^2+2( yz-xy-xt-zt), \\ \\
c:=z^2+y^2-x^2-t^2+2(xy+yt+xz-2zt), \ \ d:=x^2+y^2+z^2+t^2 ,\ \text{for some}\ x,y,z\ \text{and}\  t\in\mathbb Z.
\end{array}
\end{equation}
}

 \n Equation (\ref{maideophantineeq}) is at the heart of our main result so let us remind the reader of a few of the facts about it. For a primitive solution $[a,b,c]$ of (\ref{maideophantineeq}) we have $a,b,c\equiv \pm 1$ (mod 6).
 If we are interested in integer solutions of (\ref{maideophantineeq}) then, of course, one can permute
$[a,b,c]$ in six different ways, and change the signs, to
obtain a total of $48$ distinct ordered triples, if $a$, $b$ and $c$ are all different. If $d$ is odd, it is easy to see that $a$, $b$ and $c$ must all be  odd numbers too.  The number of all primitive solutions of (\ref{maideophantineeq}) is known (\cite{ch}, \cite{ejirt}) and it can be calculated with the formulae

\begin{equation}\label{numberofrepr}
\pi\epsilon(d):=\frac{\Lambda(d)+24\Gamma_2(d)}{48},
\end{equation}
\n where

\begin{equation}\label{numberofrep2x2plusy2}
\Gamma_2(d)=\begin{cases} 0 \ \text{ if $d$ is divisible by a prime
factor of the form $8s+5$ or $8s+7$}, \ s\ge 0,\\ \\

1 \ \text{ if $d$ is $3$}\\ \\

2^{k}\ \
\begin{cases}\text{ where $k$ is the number of
distinct prime factors of $d$ }   \\ \\
\text {of the form $8s+1$ or  $8s+3$}\  (s> 0),\end{cases}
\end{cases}
\end{equation}

\begin{equation}\label{2007hirschhorn}
\Lambda(d):=8d\underset{p|d, p\
prime}{\prod}\left(1-\frac{(\frac{-3}{p})}{p}\right),
\end{equation}

\n and $(\frac{-3}{p})$ is the Legendre symbol. We remind the reader that, if $p$ is an odd prime then

\begin{equation}\label{legendresymbol}
\left(\frac{-3}{p}\right)=\begin{cases}0\ \  {\rm if }\  p=3\\ \\
1\ \  {\rm if }\  p\equiv 1\ {\rm or\ } 7 \ {\rm (mod \ 12)} \\ \\
-1\ \  {\rm if }\  p\equiv 5\ {\rm or\ } 11 \ {\rm (mod \ 12)}
 \ \end{cases}.
\end{equation}

For some values of $d$ we may have primitive solutions of (\ref{maideophantineeq}) satisfying either $a=b$ or $b=c$.
This extra assumption leads to a Diophantine equation of the form $2a^2+c^2=3d^2$ and in  \cite{ejirt}
we have characterized the (primitive) solutions $[a,c,d]\in \mathbb N^3$ satisfying such an equation.
This was done in a manner similar to the way primitive Pythagorean triples are usually described
by a parametrization over two natural numbers. For the convenience of the reader we
include this characterization here and point out that the number of such solutions depends upon
the prime factors of $d$ which are of the form $8s+1$ or $8s+3$ (see \cite{dcox} and \cite{ejiposf}).

\begin{theorem}\label{characterizationofx=y}
Suppose that  $k$ and $\ell$ are positive integers with $k$ odd and $\gcd(k,\ell)=1$. Then $a$, $c$ and $d$ given by

\begin{equation}\label{parx=y}
d=2\ell^2+k^2 \ with \ \begin{cases} a=|2\ell^2+ 2k\ell-k^2|, \ c=|k^2+
4k\ell-2\ell^2|,\ if \ k\not\equiv \ell \ {\rm
(mod\ 3)}\\ \\
a=|2\ell^2-2k\ell-k^2|, \ c=|k^2-4k\ell-2\ell^2|,\ if \ k\not\equiv -\ell\  {\rm
(mod\ 3)}
\end{cases}
\end{equation}

\n constitute a positive primitive solution of $2a^2+c^2=3d^2$. Conversely, with the exception of the trivial solution $a=c=d=1$,
every positive primitive solution for $2a^2+c^2=3d^2$ appears in
the way described above for some $l$ and $k$.
\end{theorem}
\n The number of decompositions of $d$, as in this characterization, is well known and this number is basically included in the calculation (\ref{numberofrep2x2plusy2})
of $\Gamma_2(d)$.
Theorem~\ref{characterizationofx=y} gives us another way to construct examples of (\ref{maideophantineeq}) that are good for testing various conjectures.

In this paper we are interested in finding a simple way of constructing three points in $\mathbb Z^4$, say $A=(a_1,a_2,a_3,a_4)$,  $B=(b_1,b_2,b_3,b_4)$,
 and $C(c_1,c_2,c_3,c_4)$ with $a_i,b_i,c_i\in \mathbb Z$, $i=1,2,3,4$, such that

\begin{equation}\label{first}
\sum_{i=1}^4 (a_i-c_i)^2=\sum_{i=1}^4  (b_i-c_i)^2=\sum_{i=1}^4  (a_i-b_i)^2=D>0.
\end{equation}

\n We will disregard integer translations, so we will assume that $C$ is the origin ($C=O$). So, (\ref{first}) can be simplified to

\begin{equation}\label{first2}
\sum_{i=1}^4 a_i^2=\sum_{i=1}^4  b_i^2=\sum_{i=1}^4  (a_i-b_i)^2=D>0.
\end{equation}

A  trivial way to get new equilateral triangles from old ones is to permute the coordinates and change their signs (a total of possibly $4!(2^4)=384$ different triangles).
Also, we will mainly be talking about {\em irreducible} triangles, in the sense that the triangle $OAB$ cannot be shrunk by an integer factor to a triangle in $\mathbb Z^4$.

Our first example of a whole family of equilateral triangles in $\mathbb Z^4$ is given by any solution of (\ref{maideophantineeq}),
taking $A(d,a,b,c)$ and $B(2d,0,0,0)$.
The triangle $OAB$ is irreducible as long as we start with a primitive solution of (\ref{maideophantineeq}) and its side-lengths are equal to $2|d|$. This example shows that we
have a totally different situation compared with the 3-dimensional case where the side-lengths are never an integer. For comparison, let us include here the parametrization
from \cite{ejiobando}.

\begin{theorem}\label{oldnewparametrization}
The sub-lattice of all points in the plane ${\cal P}_{a,b,c}:=\{[\alpha,\beta,\gamma]\in \mathbb Z|a\alpha+b\beta+c\gamma=0\}$ contains two vectors $\overrightarrow{\zeta}$ and $\overrightarrow{\eta}$ such that the triangle ${{\cal T}^{m,n}_{a,b,c}}:=\triangle OPQ$ with
$P$, $Q$ in ${\cal P}_{a,b,c}$, is equilateral if and only if for
some integers $m$, $n$

\begin{equation}\label{vectorid2}\begin{array}{c}
\overrightarrow{OP}=m\overrightarrow{\zeta}-n\overrightarrow{\eta},\
\
\overrightarrow{OQ}=n\overrightarrow{\zeta}+(m-n)\overrightarrow{\eta},
\ { \rm with} \\ \\  \overrightarrow{\zeta}=(\zeta_1,\zeta_1,\zeta_2),
\overrightarrow{\varsigma}=(\varsigma_1,\varsigma_2,\varsigma_3),
\overrightarrow{\eta}=\frac{\overrightarrow{\zeta}+\overrightarrow{\varsigma}}{2},
\end{array}
\end{equation}

\begin{equation}\label{paramtwo}
\begin{array}{l}
\begin{cases}
\ds \zeta_1=-\frac{rac+dbs}{q} \\ \\
\ds \zeta_2=\frac{das-bcr}{q}\\ \\
\ds \zeta_3=r
\end{cases}
\ \ ,\ \ \ \begin{cases}
\ds \varsigma_1=\frac{3dbr-acs}{q} \\ \\
\ds \varsigma_2=-\frac{3dar+bcs}{q}\\ \\
\ds \varsigma_3=s
\end{cases}
\end{array}\ \ ,
\end{equation}
where $q=a^2+b^2$ and $r$, $s$ can be chosen so that all six
numbers in (\ref{paramtwo}) are integers. The sides-lengths of
$\triangle OPQ$ are equal to $d\sqrt{2(m^2-mn+n^2)}$. Moreover,
$r$ and $s$ can be constructed in such a way that the following
properties are also satisfied: \par

\n (i) $2q=s^2+3r^2$ and similarly
$2(b^2+c^2)=\varsigma_1^2+3\zeta_1^2$ and
$2(a^2+c^2)=\varsigma_2^2+3\zeta_2^2$\par

\n (ii) $r=r'\omega \chi$, $s=s'\omega\chi $ where
$\omega=gcd(a,b)$, $gcd(r',s')=1$ and $\chi$ is the product of the
prime factors of the form $6k-1$ of $q/\omega^2$

\n (iii) $|\overrightarrow{\zeta}|=d\sqrt{2}$,
$|\overrightarrow{\varsigma}|=d\sqrt{6}$, and
$\overrightarrow{\zeta}\cdot\overrightarrow{\varsigma} =0$

\n (iv) $s+i\sqrt{3}r=gcd(A-i\sqrt{3}B, 2q)$, in the ring $\mathbb Z[i\sqrt{3}]$, where
$A=ac$ and $B=bd$.
\end{theorem}
When it comes to the irreducible elements in $\mathbb Z[i\sqrt{3}]$ we use only the following decomposition of $4$ as
$(1+\sqrt{3}i)(1-\sqrt{3}i)$. As we can see from the above theorem, the possible values of $D$ (defined in (\ref{first}) )
for the existence of an irreducible triangle must be of the form
$2(m^2-mn+n^2)$, $m,n\in \mathbb Z$.

Let us see that in $\mathbb Z^4$ every even $D$, as in (\ref{first}), can be attained for some equilateral triangle not necessarily irreducible.
This can be accomplished by the next family given by

\begin{equation}\label{someparametrization}
A(a+b,a-b,c+d,c-d)\ and \ B(a+c,d-b,c-a,-b-d),
\end{equation}

 \n where $a$, $b$, $c$ and $d$ are arbitrary integers. The side-lengths are $\sqrt{2(a^2+b^2+c^2+d^2)}$. By Lagrange's four-square theorem
we see that for every even $D$, we can find an equilateral triangle, as in (\ref{someparametrization}), whose side-lengths are equal to $\sqrt{D}$.

 Another difference that we encounter is that in $\mathbb Z^3$, only the equilateral triangles given as in Theorem~\ref{oldnewparametrization} with $m^2-mn+n^2=k^2$, $k\in \mathbb Z$,  can be extended to a regular tetrahedron in $\mathbb Z^3$ (\cite{ejirt}). We conjecture that in $\mathbb Z^4$ every equilateral triangle can be extended to a regular tetrahedron in $\mathbb Z^4$.
 For example, the family in (\ref{someparametrization}) can be completed to a regular tetrahedron by taking
 $C(b+c,a+d,d-a,c-b)$. In this paper we will also show that every triangle in the first family given here (p. 3) can be extended to regular tetrahedrons. In fact, we will show more generally that, this is indeed the case, for the family defined by

\begin{equation}\label{parametrizationone}
\begin{array}{l}
A((m-2n)d,ma,mb,mc)\  \text{and}\ \\ \\
B((2m-n)d,na,nb,nc), \ \text{with}\  \ m,n\in \mathbb Z, \ \gcd(m,n)=1.
\end{array}
\end{equation}

\vspace{0.2in}

\section{Preliminary Results}

\vspace{0.2in}
Our classification is driven essentially by the side-length $\sqrt{D}$ introduced in (\ref{first}).

\begin{proposition}\label{irreducible} An irreducible equilateral triangle of side-lengths $\sqrt{D}$ exits if and only if
$D=2^j(2\ell-1)$ with $j\in \{1,2\}$ and $\ell \in \mathbb N$.
\end{proposition}

\n \proof.  For necessity, from (\ref{first2}) we see that $$D=|\overrightarrow{AB}|^2=\sum_{i=1}^4 a_i^2-2<\overrightarrow{OA},\overrightarrow{OB}>+\sum_{i=1}^4 b_i^2\ \   \text{which implies}\
\ \ <\overrightarrow{OA},\overrightarrow{OB}>=D/2.$$

\n So,  $\sum_{i=1}^4a_ib_i=D/2$ which proves that $D$ must be even, and so  $D=2L$ with $L\in \mathbb N$.
If $L$ is odd, then we are done, so let us assume that $L$ is a multiple of $2$.  Then $D=2L$ is a multiple of $4$. Then
$$a_1^2+a_2^2+a_3^2+a_4^2=b_1^2+b_2^2+b_3^2+b_4^2\equiv 0 \ \ (mod\ 4)$$
\n implies that either all $a_i$'s are even or all are odd, since every perfect square is equal to $0$ or $1$ (mod $4$). Similarly the $b_i$'s are all even or all odd.
 If the triangle $OAB$ is  irreducible we may assume that the $a_i$'s are all odd numbers, say $a_i=2a_i'+1$, $a_i'\in \mathbb Z$, $i=1,2,3,4$.
As a result, we have

$$D=2L=a_1^2+a_2^2+a_3^2+a_4^2=\sum_{i=1}^4(2a_i'+1)^2=
4[1+\sum_{i=1}^4a_i'(a_i'+1)].$$

\n Hence $\frac{L}{2}$ must be odd.

For sufficiency, let us use Lagrange's four square theorem and write  $2^{j-1}(2\ell-1)=a^2+b^2+c^2+d^2$, $a$, $b$, $c$ and $d$ in $\mathbb Z$. Since  $j\in \{1,2\}$, not all odd $a$, $b$, $c$ and $d$ are odd. In fact, using the Theorem~1 from \cite{ch}, (page 9), we may assume that $\gcd(a,b,c,d)=1$.  Using (\ref{someparametrization}) appropriately, we obtain an irreducible equilateral triangle $T$, whose side lengths equal $\sqrt{2^{j}(2\ell-1)}$.\eproof
\vspace{0.1in}

\n {\bf Remark:} {\em  One may want to use orthogonal matrices having rational
coordinates (see \cite{sharipov1}) in order to construct new equilateral triangles from old ones. The group of symmetries of  the four-dimensional space
which leaves the lattice $\mathbb Z^4$ invariant consists only of those transformations that permutes the variables
or change their signs. One of the simplest orthogonal matrix having rational entries is

$$O_2=\frac{1}{4}\left(
       \begin{array}{cccc}
         1 & 1 & 1 & 1 \\
         1 & 1 & -1 & -1 \\
         1 & -1 & 1 & -1 \\
         1 & -1 & -1 & 1 \\
       \end{array}
     \right),
$$
\n and  quaternion theory provides a huge class of similar examples. One can use such a matrix to transform
equilateral triangles in $\mathbb Z^4$ into equilateral triangles in $\mathbb Q^4$ and then use a dilation to get integer coordinates.This results into
changing their sides by a square factor.
By Proposition~\ref{irreducible}, we know that this method is not going to cover all the equilateral triangles
starting with a given one. In fact, we are aiming in determining a way of constructing orthogonal matrices from equilateral triangles as in the case of three dimensions. }

The next observation we make is a known way to generate new triangles from given ones. This is similar to the parametrization from Theorem~\ref{oldnewparametrization}.

\begin{proposition}\label{parammnminimality} (i) Given an equilateral triangle $OPQ$ in $\mathbb Z^4$ and $m$, $n$ two integers, then the triangle
$OP'Q'\in \mathbb Z^4$ defined by
\begin{equation}\label{parammn}
\overrightarrow{OP'}=m\overrightarrow{OP}-n\overrightarrow{OQ}, \ \ \ \overrightarrow{OQ'}=(m-n)\overrightarrow{OQ}+n\overrightarrow{OP},
\end{equation}
\n is also equilateral and its side-lengths are equal to $|\overrightarrow{PQ}|\sqrt{m^2-mn+n^2}$.

(ii) There exists an equilateral triangle $ORS\in \mathbb Z^4$  in the same plane as $\triangle OPQ$ such that all equilateral triangles in this plane (having $O$ as one of the vertices) are generated  as
in (\ref{parammn}) from $\triangle  ORS$.

(iii) If $\triangle OP'Q'$ is irreducible then $\gcd(m,n)=1$. Conversely, if $\triangle OPQ$ is minimal as in (ii), then $\triangle OP'Q'$ is  irreducible if
$\gcd(m,n)=1$.

(iv) $\triangle OPQ$ is minimal, as in part (ii), if it is irreducible and if $|\overrightarrow{PQ}|^2$ is not divisible by $3$ or by a prime of the form $6k+1$.
\end{proposition}

\n \proof.\  (i) We just calculate

$$|\overrightarrow{OP'}|^2=m^2|\overrightarrow{OP}|^2-2mn<\overrightarrow{OP},\overrightarrow{OQ}>+n^2|\overrightarrow{OQ}|^2=|\overrightarrow{PQ}|^2(m^2-mn+n^2), $$

$$|\overrightarrow{OQ'}|^2=(m-n)^2|\overrightarrow{OP}|^2+2(m-n)n<\overrightarrow{OP},\overrightarrow{OQ}>+n^2|\overrightarrow{OQ}|^2=|\overrightarrow{PQ}|^2(m^2-mn+n^2)\ and$$

$$|\overrightarrow{P'Q'}|^2=|(m-n)\overrightarrow{OP}-m\overrightarrow{OQ}|^2=|\overrightarrow{PQ}|^2(m^2-mn+n^2).$$

(ii) To show the existence of a minimal triangle,  we choose an equilateral triangle with integer coordinates $\triangle ORS$, in the same plane, which has the smallest side-lengths.
If this triangle doesn't cover all the triangles that are in the same plane, there exists one with integer coordinates, say $\triangle OAB$,  as in Figure~1, whose vertices are not part of the tessellation generated by $\triangle ORS$.
Then,  we can see that the position of $A$ and $B$ within the tessellation generated by $ORS$, provide the existence of an equilateral triangle with smaller side-lengths and having integer coordinates (in Figure 1, it is denoted $\triangle ODC$). This contradicts the assumption about the smallest side lengths.

\begin{center}
$\underset{Figure\ 1: \ Two\ distinct \ tessellations\
}{\epsfig{file=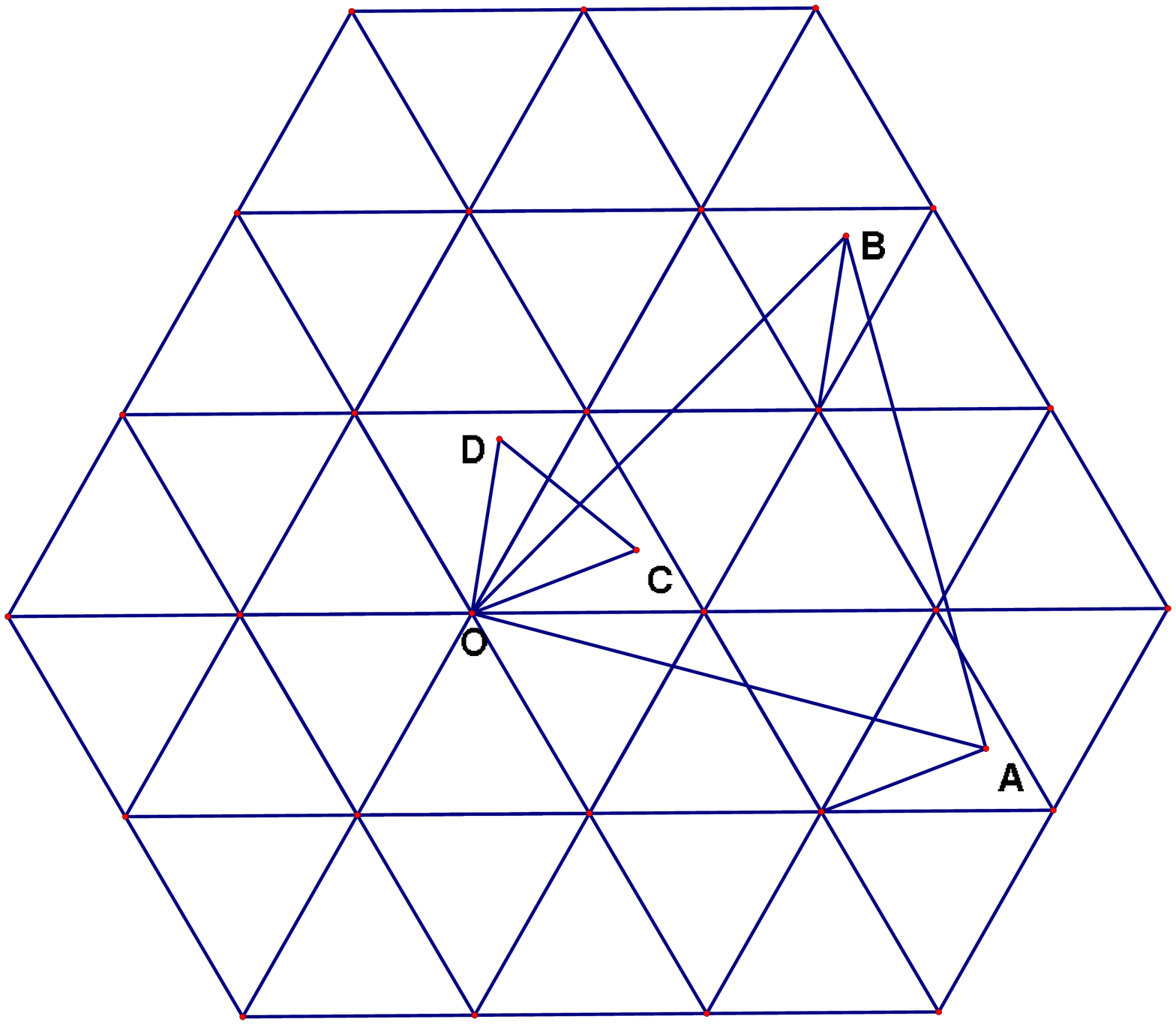,height=2in,width=2.3in}}$
\end{center}

(iii) For the first statement in part (ii), it is clear that if  $\gcd(m,n)=p$ then the triangle $OP'Q'$ can be shrunken by a factor of $p$.
If $\gcd(m,n)=1$ and $OP'Q'$ is not irreducible, then for some positive integer $t$ and an equilateral triangle $OP''Q''\in \mathbb Z^4$ we have for some $u,v\in \mathbb Z$
$$\overrightarrow{OP'}=m\overrightarrow{OP}-n\overrightarrow{OQ}=t\overrightarrow{OP'}=t(u\overrightarrow{OP}-v\overrightarrow{OQ})$$
because of the minimality of $OPQ$. We conclude that $m=tu$ and $n=tv$ which contradicts the assumption that $\gcd(m,n)=1$.

(iv) If $\triangle OPQ$ is not minimal then by (a) and (b) we see that $|\overrightarrow{PQ}|^2=(m^2-mn+n^2)\ell$ for some integers $m$, $n$ and $\ell=|\overrightarrow{RS}|^2$
where  $\triangle ORS$ is minimal. By our assumption on $|\overrightarrow{PQ}|^2$, this forces either $m^2-mn+n^2=1$ in which case $\triangle  OPQ$ is a minimal triangle, or some prime $p$ of the form $6k-1$ divides $m^2-mn+n^2$. In this case, $p$ divides $m$ and $n$ which contradicts the fact that $\triangle  OPQ$ is irreducible.\eproof

\begin{proposition}\label{aeverypoint} Suppose $A=(x_1,x_2,x_3,x_4)\in\mathbb Z^4$, different of the origin,  is a point such that the number $D=x_1^2+x_2^2+x_3^2+x_4^2$ is even. Then
there exists $B\in\mathbb Z^4$ such that $OAB$ is equilateral. Moreover, there exists $C\in\mathbb Z^4$ such that $OABC$ is a regular tetrahedron.
\end{proposition}

\n\proof.  It is clear that since $D$ is even, the number of odd numbers in the list $[x_1,x_2,x_3,x_4]$ is even. This means,
we can match them in two pairs such that the numbers in each pair are of the same parity. So,  without loss of generality, we may assume that
$x_1$ and $x_2$ have the same parity and that the same is true for $x_3$ and $x_4$. Hence the systems

$$\begin{cases}a+b=x_1\\ \\
a-b=x_2
\end{cases},\ \
\begin{cases}c+d=x_3\\ \\
c-d=x_4
\end{cases}
$$
\n have solutions in integers for $a$, $b$, $c$, and $d$. Then the statement follows from the parametrization given in (\ref{someparametrization}) and from the observation that follows (page 4). \eproof

The last result of this section is supporting the following conjecture for which we have strong numerical evidence.

\vspace{0.1in}

\n {\bf  Conjecture I:} {\it Every equilateral triangle in $\mathbb Z^4$ is a face of a regular
tetrahedron in $\mathbb Z^4$.}

\vspace{0.1in}

\begin{proposition}\label{extendthe} Conjecture I is true for the family given by (\ref{parametrizationone}).
\end{proposition}

\proof.  Let us assume that $A$ and $B$ are given by (\ref{parametrizationone}). We may try to take $R((m-n)d,x,y,z)$ such that

\begin{equation}\label{tetrahcompl}
ax+by+cz=(m+n)d^2\ \ \text{and}\ \ \ x^2+y^2+z^2=d^2(3m^2-2mn+3n^2).
\end{equation}

\n If we can find integers satisfying (\ref{tetrahcompl}), then $$|\overrightarrow{OR}|^2=x^2+y^2+z^2+(m-n)^2d^2=4d^2(m^2-mn+n^2)=|\overrightarrow{AB}|^2, $$
 $$<\overrightarrow{OR},\overrightarrow{OA}> =m(ax+by+cz)+(m-n)(m-2n)d^2 =2d^2(m^2-mn+n^2), \ \text {and}  $$

 $$<\overrightarrow{OR},\overrightarrow{OB}> =n(ax+by+cz)+(m-n)(2m-n)d^2 =2d^2(m^2-mn+n^2).$$
 This ensures that $OABR$ is a regular tetrahedron.

\n Substituting  $x=u+a\frac{m+n}{3}$, $y=v+b\frac{m+n}{3}$ and $z=w+c\frac{m+n}{3}$ where $au+bv+cw=0$, then the first equality in (\ref{tetrahcompl}) is satisfied. The second equation in (\ref{tetrahcompl})  is equivalent to $$u^2+v^2+w^2=\frac{8d^2(m^2-mn+n^2)}{3}.$$

\n We are going to use the vectors defined by the Theorem~\ref{oldnewparametrization} and take

\begin{equation}\label{definitionofuvw}
(u,v,w)=\frac{2m-4n}{3}\overrightarrow{\zeta}+\frac{2m+2n}{3}\overrightarrow{\eta},
\end{equation}

\n then

$$u^2+v^2+w^2=2d^2\left[\left(\frac{2m-4n}{3}\right)^2+\left(\frac{2m-4n}{3}\right)\left(\frac{2m+2n}{3}\right)+\left(\frac{2m+2n}{3}\right)^2\right] \Rightarrow $$

$$u^2+v^2+w^2=\frac{8d^2(m^2-mn+n^2)}{3}\ \text{and}\ au+bv+cw=0.$$

\n So, we need to determine when $x=u+a\frac{m+n}{3}$, $y=v+b\frac{m+n}{3}$ and $z=w+c\frac{m+n}{3}$ are integers if $(u,v,w)$ are defined
by (\ref{definitionofuvw}). Because of symmetry, it is enough to look at the last component; $z$ is equal to

$$z=\frac{(2m-4n)r}{3}+\frac{(2m+2n)}{3}\frac{(r+s)}{2}+c\frac{m+n}{3}=m\frac{3r+s+c}{3}-n\frac{3r-s-c}{3},$$

or

$$z=(m-n)r+(m+n)\frac{s+c}{3}.$$

Since $2q\equiv s^2\ (mod \ 3)$ and $q\equiv -c^2\ (mod \ 3)$ then $c\equiv \pm s \ (mod \ 3)$. We may choose the vector
$\overrightarrow{\varsigma}$ (changing the sign if necessary) so that $s+c\equiv 0 \ (mod \ 3)$. In fact, let us observe that we can accomplish even more by a change of sign, that
two of the coordinates  $x$, $y$  or $z$ are integers. The third coordinate must also be  an integer since $x^2+y^2+z^2=d^2(3m^2-2mn+3n^2)$. \eproof

Let us observe that the two parameterizations (\ref{someparametrization}) and (\ref{parametrizationone}) cover different classes of equilateral triangles. As in Proposition~\ref{aeverypoint}, we  can complete the points $O$ and $A=[11,5,7,17]$ in two different ways to a regular tetrahedron in $\mathbb Z^4$. If we take $d=11$ and $a=5$, $b=7$ and $c=17$, with (\ref{someparametrization}) we get  the regular tetrahedron $\{O,A, [20, -8, 4, 2], [15, 3, -13, 9]\}$.  And from (\ref{parametrizationone}) and Proposition~\ref{extendthe} we get
an essentially different regular tetrahedron $\{O,A, [22,0,0,0], [11, 1,19,-1]\}$.

%%%%%%%%%%%%%%%%%%%%%%%%%%%%%%%%%%%%%%%%%%%%%%%%%%%%%%%%%%%%%%%%%%%%%%%%%%%%%%%%%%%%%%%%%%%%%%%%%%%%%%%%%%%%%%%%%%%%%%%%%%%%
\section{Our necessary condition}
%%%%%%%%%%%%%%%%%%%%%%%%%%%%%%%%%%%%%%%%%%%%%%%%%%%%%%%%%%%%%%%%%%%%%%%%%%%%%%%%%%%%%%%%%%%%%%%%%%%%%%%%%%%%%%%%%%%%%%%%%%%%

For a set of vectors $V$ in $\mathbb R^4$ we denote as usual by $V^{\bot}$ the set of vectors $x\in \mathbb R^4$ perpendicular to every vector $v$
in $V$, i.e., $$V^{\bot}=\{[x_1,x_2,x_3,x_4]\in \mathbb R^4|x_1v_1+x_2v_2+x_3v_3+x_4v_4=0\ \text{for all}\ v=[v_1,v_2,v_3,v_4],v\in V\}.$$

\begin{proposition}\label{necessarycond} Given that $\triangle OAB$ is equilateral,
with $A=[a_1,a_2,a_3,a_4]\in  \mathbb R^4$ and
$B=[b_1,b_2,b_3,b_4]\in  \mathbb R^4$, having side lengths
$\sqrt{2L}$ ($L>0$), then there exist an odd $k$ dividing $L$ and two vectors $v$ and $w$ with integer coordinates such that

$$v=(0,\alpha_1-\beta_1,\alpha_2-\beta_2,\alpha_3-\beta_3)\
\text{and} \ w=(\alpha_3-\beta_3, -\alpha_2-\beta_2,\alpha_1+\beta_1
,0)\ \ \text{with}$$

\begin{equation}\label{thetwoeq}
\alpha_1^2+\alpha_2^2+\alpha_3^2=\beta_1^2+\beta_2^2+\beta_3^2=3k^2,
\end{equation}

\n and  $\overset{\rightarrow}{OA}, \overset{\rightarrow}{OB}\in
\{v,w\}^{\bot}$.
Moreover, we can permute the coordinates of $A$ and $B$ and/or change their
signs in order to have $v$ and $w$ linearly independent.
\end{proposition}

\n \proof. \ Using Lagrange's Identity, we get

$$4L^2=(\sum_{i=1}^4 a_i^2)(\sum_{i=1}^4  b_i^2)=(\sum_{i=1}^4a_ib_i)^2+\sum_{1\le i< j\le 4} (a_ib_j-a_jb_i)^2.$$

\n We have seen in the proof of Proposition~\ref{irreducible} that   $\sum_{i=1}^4a_ib_i=L$. This implies that

\begin{equation}\label{sixsquares}
\sum_{1\le i< j\le 4} (a_ib_j-a_jb_i)^2=3L^2.
\end{equation}

If we denote by $\Delta_{ij}=(-1)^{i-j}(a_ib_j-a_jb_i)$, we observe that $\Delta_{ij}=\Delta_{ji}$ for all $i,j\in \{1,2,3,4\}$,  $i<j$.
We have a relation which is a constraint on the solutions of (\ref{sixsquares}):

\begin{equation}\label{planeeq}
\Delta_{12}\Delta_{34}-\Delta_{13}\Delta_{24}+\Delta_{14}\Delta_{23}=0.
\end{equation}

\n This also implies that we  have a relation similar to the one in (\ref{sixsquares}), but with only three perfect squares:

\begin{equation}\label{threequares}
3L^2=(\Delta_{12}\pm \Delta_{34})^2+(\Delta_{13}\mp\Delta_{24})^2+(\Delta_{14}\pm\Delta_{23})^2.
\end{equation}

One can check that

\begin{equation}\label{planeeq2}
\begin{cases}a_1(0)+a_2\Delta_{34}+a_3\Delta_{24}+a_4\Delta_{23}=0\\
b_1(0)+b_2\Delta_{34}+b_3\Delta_{24}+b_4\Delta_{23}=0\\
a_1\Delta_{23}+a_2\Delta_{13}+a_3\Delta_{12}+a_4(0)=0\\
b_1\Delta_{23}+b_2\Delta_{13}+b_3\Delta_{12}+a_4(0)=0.
\end{cases}
\end{equation}

So, to prove our statement we take $\alpha_1=\Delta_{12}+
\Delta_{34}$, $\alpha_2=-\Delta_{13}+\Delta_{24}$,
$\alpha_3=\Delta_{14}+\Delta_{23}$, $\beta_1=\Delta_{12}-
\Delta_{34}$, $\beta_2=-\Delta_{13}-\Delta_{24}$, and
$\beta_3=\Delta_{14}-\Delta_{23}$. In general, an equation as (\ref{threequares}) can be simplified by a factor $\ell$ of $L$ and we can always assume that
$k=L/\ell$ is odd. For instance, if $L$ is even, then (\ref{threequares}) can be simplified by $4$. Any such possible simplification does not change the space $\{v,w\}^{\bot}$ since the new vectors are simply scalar multiples of $v$ and $w$. Or, in the case of the family (\ref{parametrizationone}), since $L=2d^2(m^2-mn+n^2)$, we have $\ell=2d(m^2-mn+n^2)$ and $k=d$.

It is easy to check that if $\alpha_3\not = \beta_3$ then
the two vectors $v$ and $w$ are linearly independent and so the
orthogonal space  $\{v,w\}^{\bot}$ is two-dimensional. By permuting
 coordinates and/or changing their signs, we can insure that $\alpha_3-\beta_3=2\Delta_{23}\not=0$, since by (\ref{sixsquares}) and the assumption that
$k>0$, we see that not all of the $\Delta_{ij}$ can be zero.\eproof

\n {\bf Remark I:} Let us observe that the  Gram determinant of the vectors $v$ and $w$, defined as in Proposition~\ref{necessarycond}, has an interesting value:

$$\begin{array}{c}G:=\left|
    \begin{array}{cc}
      <v,v> & <v,w> \\
      <w,v> & <w,w> \\
    \end{array}
  \right|=\left|
    \begin{array}{cc}
      6k^2-2(\alpha_1\beta_1+\alpha_2\beta_2+\alpha_3\beta_3) & 2(\alpha_2\beta_1-\alpha_1\beta_2) \\
      2(\alpha_2\beta_1-\alpha_1\beta_2) &  6k^2+2(\alpha_1\beta_1+\alpha_2\beta_2-\alpha_3\beta_3) \\
    \end{array}
    \right|= \\ \\ =(6k^2-2\alpha_3\beta_3)^2-4(\alpha_1^2\beta_1+\alpha_2\beta_2)^2-4(\alpha_2\beta_1-\alpha_1\beta_2)^2=\\ \\
    (6k^2-2\alpha_3\beta_3)^2-4(\alpha_1^2+\alpha_2^2)(\beta_1^2+\beta_2^2)=36k^4-24k^2\alpha_3\beta_3+4\alpha_3^2\beta_3^2-4(3k^2-\alpha_3^2)(3k^2-\beta_3^2)=\\ \\
    12k^2(\alpha_3^2+\beta_3^2-2\alpha_3\beta_3)=12k^2(\alpha_3-\beta_3)^2=48k^2\Delta_{23}^2. \end{array}$$
This calculation has something to do with the area of the fundamental domain of the lattice of points determined by all the points of integer coordinates
in the orthogonal subspace determined by the two vectors $v$ and $w$. Indeed, if $\Delta_{23}=1$ is equal to one, then this subspace is generated by two vectors,
$v'=[(\alpha_2+\beta_2)/2,1,0,(\beta_1-\alpha_1)/2]$ and $w'=[-(\alpha_1+\beta_1)/2,0,1,(\beta_2-\alpha_2)/2]$, assuming that these two vectors have integer coordinates also.
Hence the area of the parallelogram determined by these two vectors is $A=|v'||w|\sin(\psi)$ where $\psi$ is the angle between $v'$ and $w'$. Since $\cos(\psi)=\frac{<v',w'>}{|v'||w'|}$ then $A=\sqrt{|v'|^2|w'|^2-<v',w'>^2}=\sqrt{G/16}=k\sqrt{3}$. At the end of the paper we will exemplify more about this and its connection with the Ehrhart polynomial associated to those triangles. This is very similar to the 3-dimensional situation (see \cite{ejiep}).

\vspace{0.1in}

{\bf Remark II:} Let us observe that in example (\ref{someparametrization}) we have $\Delta_{12}=a^2+ac-bc-ad-db+b^2$, $\Delta_{34}=c^2-ac+ad+bc+db+d^2$, $\Delta_{13}=-a^2+bc-ab-c^2-ad-cd$,
$\Delta_{14}=ac+c^2-cd+ab+b^2+db$,  $\Delta_{23}=cd+d^2-db-ac+a^2-ab$, and $\Delta_{24}=-ab-ad+b^2-cd+bc+d^2$. We also see that

$$\Delta_{12}+\Delta_{34}=-\Delta_{13}+\Delta_{24}=\Delta_{14}+\Delta_{23}=a^2+b^2+c^2+d^2=k.$$

\n We observe that these calculations give us a candidate for a parametrization of all of the solutions of the Diophantine equation (\ref{maideophantineeq}) which was mentioned in (\ref{parofmain}).

\vspace{0.1in}

\vspace{0.1in}

\section{Our necessary condition is also sufficient}

\vspace{0.1in}

\begin{theorem}\label{maintheorem} Given $k$ odd, and two different representations of (\ref{maideophantineeq}), i.e.,
$$3k^2=a^2+b^2+c^2=a'^2+b'^2+c'^2,\ \text{with}\ \gcd(a,b,c,a',b',c')=1,\ c'>c, $$

\n   then if we set $\Delta_{12}=\frac{a'-a}{2}$, $\Delta_{34}=\frac{a+a'}{2}$, $\Delta_{13}=-\frac{b'-b}{2} $, $\Delta_{24}=\frac{b+b'}{2}$,
$\Delta_{14}=\frac{c+c'}{2}$, and $\Delta_{23}=\frac{c'-c}{2}$  the equations (\ref{planeeq}) and (\ref{threequares}) are satisfied.
Moreover, the two  dimensional space
$\cal S$ of all vectors $[u,v,w,t]\in \mathbb Z^4$, such that

\begin{equation}\label{planeeqsc}
\begin{cases}(0)u+\Delta_{34}v+\Delta_{24}w+\Delta_{23}t=0\\
 \Delta_{23}u+\Delta_{13}v+\Delta_{12}w+(0)t=0
\end{cases}
\end{equation}

\n contains a family of equilateral triangles as in Proposition~\ref{parammnminimality}. Conversely, every equilateral triangle in $\mathbb Z^4$ can be obtained in this way by Proposition~\ref{necessarycond}.
\end{theorem}

\n \proof.\  Clearly by construction we see that (\ref{threequares}) is true. Since $a$, $a'$, $b$, $b'$, $c$ and $c'$, are all odd integers, then all numbers defined above are integers.  We notice that $$\Delta_{12}\Delta_{34}-\Delta_{13}\Delta_{24}+\Delta_{14}\Delta_{23}=(1/4)[a'^2-a^2-(b^2-b'^2)+(c'^2-c^2)]
=0,\ \text{and}$$

$$\sum_{i,j}\Delta_{ij}^2=\frac{1}{4}[(a-a')^2+(a+a')^2+(b-b')^2+(b+b')^2+(c-c')^2+(c+c')^2]=3k^2.$$

\n Let us introduce $\Delta_1:=\gcd( \Delta_{34},\Delta_{24},\Delta_{23})$ and $\Delta_2:=\gcd(\Delta_{23},\Delta_{13},\Delta_{12})$.
Clearly $\gcd(\Delta_1,\Delta_2)=1$. The equations of the hyper-planes in (\ref{planeeqsc}) do not depend on $\Delta_1$ and $\Delta_2$. We see that the assumption
$\Delta_{23}>0$ insures that the equations (\ref{planeeqsc}) define a two dimensional space in $\mathbb R^4$.  In order to prove our claim, it suffices to show that an equilateral triangle with rational coordinates exists in
$\cal S$. We can solve the equations (\ref{planeeqsc}) for rational values of $t$ and $u$:

$$t=-\frac{\Delta_{34}v+\Delta_{24}w}{\Delta_{23}}, \ u= -\frac{\Delta_{13}v+\Delta_{12}w}{\Delta_{23}},\ \text{with}\ v,w\in \mathbb Z.$$

\n If we denote a generic point $P\in\mathbb R^4$ in the plane (\ref{planeeqsc}), in terms of $v$ and $w$, i.e.,

$$P(v,w)=[-\frac{\Delta_{13}v+\Delta_{12}w}{\Delta_{23}}, v,w, -\frac{\Delta_{34}v+\Delta_{24}w}{\Delta_{23}}],$$

\n we observe that for two pairs of the parameters, $(v,w)$ and $(v',w')$, we obtain $\Delta_{23}'=wv'-vw'$,
$$\Delta_{12}'=v(-\frac{\Delta_{13}v'+\Delta_{12}w'}{\Delta_{23}})+\frac{\Delta_{13}v+\Delta_{12}w}{\Delta_{23}}v'=\Delta_{12}\Delta_{23}'/\Delta_{23},$$
\n and similarly $\Delta_{13}'=\Delta_{13}\Delta_{23}'/\Delta_{23}$, etc.

We will show that a non-zero equilateral triangle of the form $OP(v,w)P(v',w')$ exists for some integer values of  $v$, $w$, $v'$ and $w'$. From the proof of Proposition~\ref{necessarycond}, we see that it is sufficient to have $OP(v,w)^2=OP(v',w')^2=2L$ and $$<OP(v,w),OP(v',w')>=L$$ for some $r>0$. By Proposition~\ref{necessarycond}, we may take  $L=k\ell$ for some $\ell \in \mathbb Z$, $\ell>0$, since the equations (\ref{thetwoeq}) may be simplified by some integer.

By Lagrange's identity, calculations similar to those  in the derivation of (\ref{sixsquares}), show that if
$|P(v,w)|^2=|P(v',w')|^2=2k\ell$ then

$$(<P(v,w),P(v',w')>)^2=4k^2\ell^2-(\sum_{i,j}\Delta_{ij}^2)\left(\frac{\Delta_{23}'}{\Delta_{23}}\right)^2=4k^2\ell^2-3k^2\left(\frac{\Delta_{23}'}{\Delta_{23}}\right)^2.$$

Hence, to get $<OP(v,w),OP(v',w')>=k\ell$, it is enough to have $\Delta_{23}'=\ell \Delta_{23}$ (provided that we have already shown that $|P(v,w)|^2=|P(v',w')|^2=2k\ell$).
In other words, an equilateral triangle exists if there exist integer pairs $(v,w)$ and $(v',w')$ such that

\begin{equation}\label{equivfrmulation}
|OP(v,w)|^2=|OP(v',w')|^2=2k\frac{\Delta_{23}'}{\Delta_{23}}.
\end{equation}

So, we will concentrate on this goal for the rest of the proof.

\n For infinitely many values of $v$ and $w$ we get integer values for $t$ and $u$, but for calculation purposes we will simply work with rational values and consider the quadratic form

$$QF(v,w):=|OP(u,v)|^2=\left(\frac{\Delta_{34}v+\Delta_{24}w}{\Delta_{23}}\right)^2+v^2+w^2+\left(\frac{\Delta_{13}v+\Delta_{12}w}{\Delta_{23}}\right)^2,$$

or

$$QF(v,w)=\frac{(\Delta_{34}^2+\Delta_{13}^2+\Delta_{23}^2)v^2+2(\Delta_{34}\Delta_{24}+\Delta_{13}\Delta_{12})vw
+(\Delta_{24}^2+\Delta_{12}^2+\Delta_{23}^2)w^2}{\Delta_{23}^2}.$$

\n We use Lagrange's identity again, which we will write in the form

$$(\alpha^2+\beta^2+\gamma^2)(\zeta^2+\eta^2+\theta^2)=(\alpha\zeta-\beta\eta+\gamma\theta)^2
+(\alpha\eta+\beta\zeta)^2+(\alpha\theta-\gamma\zeta)^2+(\beta\theta+\gamma\eta)^2.$$

\n Then the determinant of the form $QF(v,w)$ (excluding its denominator) is equal to

$$\begin{array}{c}
-\Delta/4=(\Delta_{34}^2+\Delta_{13}^2+\Delta_{23}^2)(\Delta_{12}^2+\Delta_{24}^2+\Delta_{23}^2)-(\Delta_{34}\Delta_{24}+\Delta_{13}\Delta_{12})^2=\\ \\
(\Delta_{12}\Delta_{34}-\Delta_{13}\Delta_{24}+\Delta_{23}^2)^2+(\Delta_{12}\Delta_{23}-\Delta_{34}\Delta_{23})^2+(\Delta_{13}\Delta_{23}+\Delta_{24}\Delta_{23})^2\Rightarrow\\ \\
-\Delta/4 =\Delta_{23}^2\left[(\Delta_{23}-\Delta_{14})^2+(\Delta_{12}-\Delta_{34})^2+(\Delta_{13}+\Delta_{24})^2 \right]=\Delta_{23}^2(3k^2)\\ \\  \Rightarrow \Delta=-3(2k\Delta_{23})^2.
\end{array}$$

\n This implies that for $v_0= -(\Delta_{34}\Delta_{24}+\Delta_{13}\Delta_{12})$ and $w_0=\Delta_{13}^2+\Delta_{23}^2+\Delta_{34}^2$ we get

$$QF(v_0,w_0)=3k^2(\Delta_{34}^2+\Delta_{13}^2+\Delta_{23}^2),$$

\n and

\begin{equation}\label{qf}
QF(v,w)=\frac{(w_0v-v_0w)^2+3k^2w^2\Delta_{23}^2}{\Delta_{23}^2w_0}.
\end{equation}

\n So, by (\ref{equivfrmulation}) and (\ref{qf}) we need to have solutions $(v,w)$ and $(v',w')$ of

 \begin{equation}\label{suffices}
(w_0v-v_0w)^2+3k^2w^2\Delta_{23}^2 =(w_0v'-v_0w')^2+3k^2w'^2\Delta_{23}^2 =2k\Delta_{23}w_0(wv'-vw').
 \end{equation}

 \n In general, if the Diophantine quadratic equation $x^2+3y^2=n$ has a solution $(x,y)$, then it has other solutions other than the trivial ones (changing the signs of $x$ and $y$). For instance, we have
 $$\left(\frac{x+3y}{2}\right)^2+3\left(-\frac{x-y}{2}\right)^2=x^2+3y^2=n.$$

\n  Given $(v,w)$, let us assume that $(v',w')$ gives exactly this other solution ($x=w_0v-v_0w$, $y=kw\Delta_{23}$) when substituted in (\ref{suffices}):

$$\ds
\begin{cases}\ds w_0v'-v_0w'=\frac{w_0v-v_0w+3kw\Delta_{23}}{2}\\ \\
\ds k\Delta_{23}w'=-\frac{w_0v-v_0w-kw\Delta_{23}}{2}.
\end{cases}
$$

\n This allows us to solve for $v'$ and $w'$ in order to calculate $\Delta_{23}'$:

$$v'=\frac{v_0^2w+(k\Delta_{23}-v_0)w_0v+3k^2\Delta_{23}^2w}{2k\Delta_{23}w_0},\ w'=\frac{v_0w-w_0v}{2k\Delta_{23}}+\frac{w}{2},$$

$$\Delta_{23}'=wv'-vw'=\frac{(w_0v-v_0w)^2+3k^2w^2\Delta_{23}^2}{2k\Delta_{23}w_0}.$$

\n This means that (\ref{suffices}) is automatically satisfied with this choice of $(v',w')$. In fact, we have shown that, for every integer values of
$(v,w)$, such that

\begin{equation}\label{criticaleq}
(w_0v-v_0w)^2+3k^2w^2\Delta_{23}^2 =(2kw_0\ell) \Delta_{23}^2
\end{equation}

\n taking $(v',w')$ as above, we automatically get an equilateral triangle, $OP(v,w)P(v',w')\in \mathbb Q^4$, which may not have integer coordinates.
In order for (\ref{criticaleq}) to have rational solutions, in $v$ and $w$, we need to have $2kw_0\ell$ a positive integer which is a product of
primes of the form $6s+1$ or equal to $3$, and all the other primes, in its prime factorization, have even exponents. Clearly, there exists the smallest positive integer $\ell$ with this property. \eproof

Next, let us look at some examples. The analysis is quite simple if we start with a prime $k$, say $k=11$, which is not a prime of the form $6s+1$.
We can take all sorts of possible solutions for (\ref{maideophantineeq}) ($d=k=11$) as pointed out in the Introduction, so let's say we choose $a=1$, $b=1$, $c=-19$, $a'=5$, $b'=7$, $c'=-17$. Then,  $\Delta_{12}=\frac{a'-a}{2}=2$, $\Delta_{34}=\frac{a+a'}{2}=3$, $\Delta_{13}=-\frac{b'-b}{2}=-3 $, $\Delta_{24}=\frac{b+b'}{2}=4$,
$\Delta_{14}=\frac{c+c'}{2}=-18$, and $\Delta_{23}=\frac{c'-c}{2}=1$. We would like to determine the minimal triangles defined by this data. This triangle is contained in the plane defined the by the system of equations:

\begin{equation}\label{planeeqscex1}
\begin{cases}3v+4w+t=0\\
u-3v+2w=0.
\end{cases}
\end{equation}

This gives the quadratic form
$$\begin{array}{c}QF(v,w)=(3v+4w)^2+v^2+w^2+(3v-2w)^2=19v^2+12uv+21w^2=\\ \\ \ds =\frac{(19v+6w)^2+3(11)^2w^2}{19}
\end{array}$$
\n which implies that the equation (\ref{criticaleq}) becomes $(19v+6w)^2+3(11)^2w^2=2(11)(19)\ell$ and so the smallest $\ell$ that insures solutions,
at least for rational numbers $v$ and $w$, is $\ell=22$. Since $19=4^2+3(1^2)$ we see that $v=4$ and $w=2$ work out as integer solutions of (\ref{criticaleq}).
It turns out that $v'$ and $w'$ are also integers and so we get $P=[8, 4, 2, -20]$ and $P'=[21, 5, -3, -3]$ that gives us a minimal triangle, $\triangle_1=OPP'$ (see Proposition~\ref{parammnminimality}), of side-lengths equal to
$22$, and contained in the plane (\ref{planeeqscex1}).  It turns out, according to our calculations, that this triangle is just one of the 40 essentially different triangles of side-lengths $22$ (see Table~1 below for other examples). As we have pointed out in the Remark I after Proposition~\ref{necessarycond}, the fundamental domain of the lattice
of integer points in the plane containing this triangle has an area of $11\sqrt{3}$. Therefore, its Ehrahrt polynomial is equal to

$$P(\triangle_1,t)=11t^2+2t+1,\ t\in \mathbb N,$$

\n which implies that there are exactly 10 points of the lattice in its interior.

Let us take an example when $k$ is a composite number, say $k=15$. We will take $a=1$, $b=7$, $c=25$, $a'=3$, $b'=15$, and $c'=21$. Then
 $\Delta_{12}=\frac{a'-a}{2}=1$, $\Delta_{34}=\frac{a+a'}{2}=2$, $\Delta_{13}=-\frac{b'-b}{2}=-4 $, $\Delta_{24}=\frac{b+b'}{2}=11$,
$\Delta_{14}=\frac{c+c'}{2}=23$, and $\Delta_{23}=\frac{c'-c}{2}=2$. A minimal triangle for this choice is contained in the plane:

\begin{equation}\label{planeeqscex1}
\begin{cases}2v+11w+2t=0\\
2u-4v+w=0.
\end{cases}
\end{equation}

 We see that $w=2w'$ and so we have a quadratic form
$$\begin{array}{c}QF(v,w')=(v+11w')^2+v^2+4w'^2+(2v-w')^2=6v^2+18vw'+126w'^2=\\ \\ \ds =3\frac{4v^2+12vw'+84w'^23}{2}=
3\frac{(2v+3w')^2+3(25)w'^2}{2}.
\end{array}$$
 Since $w_0=24$, we need to solve the Diophantine equation $ (2v+3w')^2+3(25)w'^2=20\ell $. The smallest $\ell $ for which we have solutions is $\ell=5$.
We can see that $v=w'=1$ is a possible solution which gives $P=[ 1, 1, 2, -12]$. This gives $P'=[10,5,0,-5]$, a minimal triangle $\triangle_2=OPP'$ of side lengths equal to
$150$. We observe that (\ref{planeeqscex1}) is equivalent to

$$\begin{cases}v+11w'+t=0\\
u-2v+w'=0.
\end{cases}$$

\n which shows that the lattice of integer points satisfying this system is generated by the vectors $[2,1,0,-1]$ and $[-1,0,2,-11]$.
In this case the area of the fundamental domain is $A=15\sqrt{3}$, and so its Ehrhart polynomial is

$$P(\triangle_2,t)=25t^2+7t+1,\ t\in \mathbb N.$$

We implemented the characterization from Theorem~\ref{maintheorem} numerically and calculated all of the \underline{minimum}
equilateral triangles of a given side, up to translations and symmetries of the smallest box, $[0,m_1]\times [0,m_2]\times [0,m_3]\times [0,m_4] $,
containing the equilateral triangle. The results for the side-lengths $\sqrt{2L}$ with $L=1,2,...,21$ and the corresponding values of $k$ are included in the following table.
We use a representation of the triangles which has the origin as one the vertices and the other two points have as many as possible negative
coordinates (the only points listed). In general the connection between $L$ and $k$ seems to be complicated but numerical evidence shows that
$k\le L\le 2k^2$.

$$\overset{Table 1}{\\  \begin{tabular}{|c|c|p{13cm}|c|}
  \hline  \hline
  L & \# &  \centerline{Minimal Triangles} & k \\
   \hline\underline{}
  1 & 1 & \{[0,0,1,1],[0,1,1,0]\}&  1   \\
   \hline
  2 & 1 & \{[0,0,0,2],[1,1,1,1]\}& 1 \\
    \hline
  3 & 1 & \{[0,1,1,2],[1,1,2,0]\}& 3 \\
      \hline
   4 & 0 &  &   \\
    \hline
    5 & 1 & \{[0,0,1,3],[1,2,2,1]\} &  5 \\
    \hline
     6 & 1 & \{[0,2,2,2],[1,1,3,-1]\} &  5 \\
    \hline
    7 & 1 & \{[0,1,2,3],[1,3,2,0]\} &  7 \\
    \hline
     8 & 0 &  &   \\
    \hline
    9 & 4 & \{[1,2,2,3],[0,4,-1,1]\},\{[0,1,1,4],[0,1,4,1]\},\{[0,0,3,3],[2,3,2,1] \},\{[0,1,1,4],[2,2,3,1]\} &  3, 9 \\
    \hline
    10 & 1 & \{[0,0,2,4],[1,3,3,1]\} &  5 \\
    \hline
    11 & 3 & \{[0,2,3,3],[2,1,4,-1]\},\{[0,2,3,3],[1,4,2,-1]\},\{[0,2,3,3],[4,1,2,1]\} &  11 \\
    \hline
     12 & 0 &  &   \\
        \hline
      13   & 1 & \{[0,1,3,4],[1,4,3,0]\} &  13 \\
    \hline
     14   & 1 & \{[1,3,3,3],[2,2,4,-2]\} &  7 \\
    \hline
     15   & 6 & \{[0,1,2,5],[1,5,0,2]\},\{[0,1,2,5],[2,4,3,1]\},\{[0,1,2,5],[4,3,1,2]\}, \break 
     \{[1,2,3,4],[1,4,-2,3]\},\{[1,2,3,4],[3,1,-2,4]\},\{[1,2,3,4],[4,3,-1,2]\}  &   15 \\
    \hline
    16 & 0 &  &   \\
        \hline 
        17 & 4 & \{[0,0,3,5],[1,4,4,1]\},\{[0,3,3,4],[2,2,5,-1]\},\{[0,3,3,4],[1,5,-2,2]\},\{[1,2,2,5],[3,3,4,0]\}&  17 \\
    \hline
    18 & 2 & \{[0,2,4,4],[1,5,-1,3]\},\{[0,2,4,4],[5,1,1,3]\}  & 9 \\    \hline 
    19   & 4 &\{[0,1,1,6],[3,3,4,2]\},\{[0,2,3,5],[2,5,3,0]\},\{[2,3,3,4] [2,4,-3,3]\},\{[2,3,3,4],[3,2,-3,4]\} &  19 \\
    \hline
    20 & 0 &  &   \\
        \hline 
   21   & 6 &\{[0,1,1,6],[3,3,4,2]\},\{[0,1,4,5],[4,4,3,1]\},\{[0,1,4,5],[1,5,4,0]\},\break 
   \{[0,1,4,5],[5,3,2,2]\},\{[1,1,2,6],[1,2,6,1]\},\{[1,3,4,4],[2,5,3,-2]\}  &   21 \\   \hline 
        \end{tabular}}$$


\begin{thebibliography}{99}
\bibitem{rceji} R.~Chandler and E.~J.~Ionascu, {\it A characterization of all equilateral triangles in $\Bbb Z^3$}, Integers, Art. A19 of Vol.  8({\bf 2008})

\bibitem{ch} S.~Cooper and M.~Hirschhorn, {\it On the number of primitive representations of integers as a sum of squares},  Ramanujan J. 13 ({\bf 2007}), pp. 7-25

\bibitem{dcox} D.~A.~Cox, {\em  Primes of the Form $x^2 + ny^2$}, Wiley-Interscience, {\bf 1989}


\bibitem{eji} E.~J.~Ionascu, {\it A parametrization of equilateral triangles having integer coordinates}, Journal of Integer Sequences, Vol. {\bf 10}, 09.6.7
(2007)

\bibitem{ejiobando} E. J. Ionascu and R. Obando, {\em Counting all cube in $\{0,1,2,...,n\}^3$},  Integers, 12A ({\bf 2012}) (John Selfridge Memorial Issue), Art. A9


\bibitem{ejirt} E. J. Ionascu, {\it A characterization of regular tetrahedra in $\Bbb Z^3$},
J. Number Theory, 129({\bf 2009}), pp. 1066-1074

\bibitem{ejiposf} E. J. Ionascu, {\it Primes of the form $\pm a^2\pm qb^2$},  ({\bf 2012}), arXiv:1207.0172v1 (update in preparation)

\bibitem{ejiep} E. J. Ionascu, {\it  Ehrhart's polynomial for equilateral triangles in $\Bbb Z^3$}, Australas. J. Combinatorics, 55 ({\bf 2013}),
pp. 189-204

\bibitem{ejirtlv} E. J. Ionascu, {\it Regular tetrahedra with integer coordinates of their vertices},
{\em Acta Math. Univ. Comenianae}, Vol. LXXX, 2({\bf 2011}), pp. 161-170

\bibitem{mc} I. G. McDonald, {\it Regular simplexes with integer vertices}, C. R. Math. Rep. Acad. Sci. Canada, Vol IX, No 4, ({\bf 1987}), pp. 189-193

\bibitem{sharipov1} R. A. Sharipov, {\it Algorithm for generating orthogonal matrices with rational elements}, ({\bf 2009}), arXiv:cs/0201007v2

\bibitem{s} I. J. Schoenberg, {\it Regular Simplices and Quadratic Forms}, J. London Math. Soc. 12 ({\bf 1937}), pp. 48-55

\end{thebibliography}
\end{document}